\documentclass[a4paper]{scrartcl}

\usepackage[utf8]{inputenc}
\usepackage[T1]{fontenc}

\usepackage{amsmath, amsthm, amssymb}
\usepackage{mathtools}
\usepackage{url}
\usepackage{array,booktabs,multirow}
\usepackage{subcaption}
\DeclareCaptionLabelFormat{subtab}{\tablename #1~\arabic{table}\alph{subtable} \textbf{--}}
\clearcaptionsetup{table}
\usepackage{enumitem}
\newlist{hypothenum}{enumerate}{3}
\setlist[hypothenum,1]{label=(\roman*)}

\theoremstyle{plain}

\newtheorem{theorem}{Theorem}[section]
\newtheorem{proposition}[theorem]{Proposition}

\newtheorem{problem}[theorem]{Problem}
\newtheorem{conjecture}[theorem]{Conjecture}

\theoremstyle{definition}

\theoremstyle{remark}

\numberwithin{equation}{section}

\newcommand{\setsuch}[2]{\left\{ #1 \; \middle| \; #2 \right\}}
\newcommand{\restr}[2]{{\left. #1 \right|}_{#2}}

\newcommand{\RR}{\mathbb{R}}
\newcommand{\CC}{\mathbb{C}}

\newcommand{\lie}{\mathfrak}

\DeclareMathOperator{\Id}{Id}

\DeclareMathOperator{\SO}{SO}

\usepackage[all,cmtip]{xy}
\newcommand{\fundef}[5]{
\entrymodifiers={+!!<0pt,\fontdimen22\textfont2>}
\xymatrix@R=3pt{\llap{$#1$\;\;} {#2} \ar@{->}[r] & {#3} \\ {#4} \ar@{|->}[r] & {#5}}
} 

\newcommand{\ie}{i.e.\ }
\newcommand{\eg}{e.g.\ }

\begin{document}

\title{Action of $w_0$ on $V^L$ for orthogonal and exceptional groups}
\author{Ilia Smilga}
              
\maketitle

\begin{abstract}
In this note, we present some results that partially answer the following question. Let $G$ be a simple real Lie group; what is the set of representations~$V$ of~$G$ in which the longest element~$w_0$ of the restricted Weyl group~$W$ acts nontrivially on the subspace~$V^L$ of~$V$ formed by vectors that are invariant by~$L$, the centralizer of a maximal split torus of~$G$? We give a conjectural answer to that question, as well as the experimental results that back this conjecture, when $G$ is either an orthogonal group (real form of $\SO_n(\mathbb{C})$ for some $n$) or an exceptional group.
\end{abstract}

\section{Introduction and motivation}

\subsection{Basic notations and statement of problem}

Let $G$ be a semisimple real Lie group, $\lie{g}$ its Lie algebra, $\lie{g}^\CC$ the complexification of~$\lie{g}$. We start by establishing the notations for some well-known objects related to~$\lie{g}$.

\begin{itemize}
\item We choose in~$\lie{g}$ a \emph{Cartan subspace}~$\lie{a}$ (an abelian subalgebra of $\lie{g}$ whose elements are diagonalizable over~$\RR$ and which is maximal for these properties).
\item We choose in~$\lie{g}^\CC$ a \emph{Cartan subalgebra}~$\lie{h}^\CC$ (an abelian subalgebra of $\lie{g}^\CC$ whose elements are diagonalizable and which is maximal for these properties) that contains~$\lie{a}$.
\item We denote $L := Z_G(\lie{a})$ the centralizer of~$\lie{a}$ in~$G$, $\lie{l}$ its Lie algebra.
\item Let $\Delta$~be the set of roots of~$\lie{g}^\CC$ in~$(\lie{h}^\CC)^*$. We shall identify $(\lie{h}^\CC)^*$ with~$\lie{h}^\CC$ via the Killing form. We call $\lie{h}_{(\RR)}$ the $\RR$-linear span of~$\Delta$; it is given by the formula $\lie{h}_{(\RR)} = \lie{a} \oplus i \lie{a}^\perp$.
\item We choose on~$\lie{h}_{(\RR)}$ a lexicographical ordering that ``puts $\lie{a}$ first'', \ie such that every vector whose orthogonal projection onto~$\lie{a}$ is positive is itself positive. We call $\Delta^+$ the set of roots in~$\Delta$ that are positive with respect to this ordering, and we let $\Pi = \{\alpha_1, \ldots, \alpha_r\}$ be the set of simple roots in~$\Delta^+$ (where $r$ is the rank of $\lie{g}^\CC$). Let~$\varpi_1, \ldots, \varpi_r$ be the corresponding fundamental weights.
\item We introduce the \emph{dominant Weyl chamber}
\[\lie{h}^+ := \setsuch{X \in \lie{h}}{\forall i = 1, \ldots, r,\quad \alpha_i(X) \geq 0},\]
and the \emph{dominant restricted Weyl chamber}
\[\lie{a}^+ := \lie{h}^+ \cap \lie{a}.\]
\item We introduce the \emph{restricted Weyl group} $W := N_G(\lie{a})/Z_G(\lie{a})$ of~$G$. Then $\lie{a}^+$ is a fundamental domain for the action of~$W$ on~$\lie{a}$. We define the \emph{longest element} of the restricted Weyl group as the unique element $w_0 \in W$ such that $w_0(\lie{a}^+) = -\lie{a}^+$.
\item For each \emph{dominant integral weight} $\lambda$ of $\lie{g}^\CC$ (\ie linear combination of the fundamental weights $\varpi_i$ with nonnegative integer coefficients), we denote by $V_\lambda$ the irreducible representation of $\lie{g}$ with highest weight $\lambda$.
\end{itemize}

Our goal is to study the action of~$W$, and more specifically of~$w_0$, on various representations~$V$ of~$G$. Note however that this action is ill-defined: indeed if we want to see the abstract element $w_0 \in W = N_G(\lie{a})/Z_G(\lie{a})$ as the projection of some concrete map $\tilde{w_0} \in N_G(\lie{a}) \in G$, then $\tilde{w}_0$ is defined only up to multiplication by an element of $Z_G(\lie{a}) = L$, whose action on~$V$ can of course be nontrivial.

This naturally suggests the idea of restricting to $L$-invariant vectors. Given a representation $V$ of~$\lie{g}$, we denote
\[V^L := \setsuch{v \in V}{\forall l \in L,\;\; l \cdot v = v}\]
the $L$-invariant subspace of~$V$: then $W$, and in particular $w_0$, has a well-defined action on~$V^L$.

Our goal is to characterize, for a given semisimple real Lie group~$G$, the representations $V$ of~$G$ for which the action of~$w_0$ on $V^L$ is nontrivial. This problem naturally splits into two subproblems (see \cite{Smi20} for a more extended discussion):
\begin{problem}
\label{nontrivial_Vl}
Given a semisimple Lie algebra~$\lie{g}$ and a dominant integral weight~$\lambda$, give a simple necessary and sufficient condition for having $V_\lambda^\lie{l} \neq 0$.
\end{problem}
\begin{problem}
\label{nontrivial_w0_action_simple}
Given a simple Lie algebra~$\lie{g}$ and a dominant integral weight~$\lambda$, assuming that $V_\lambda^\lie{l} \neq 0$, give:
\begin{hypothenum}
\item a simple necessary and sufficient condition for having $\restr{w_0}{V_\lambda^\lie{l}} = \pm \Id$;
\item a criterion to determine the actual sign.
\end{hypothenum}
\end{problem}
In~\cite{Smi20b}, we have already completely solved problem~\ref{nontrivial_Vl}. In~\cite{LFlSm}, we have solved problem~\ref{nontrivial_w0_action_simple} in the case where $\lie{g}$ is split. In this note, we shall present our recent work on this latter problem.

\subsection{Background and motivation}

These two problems arose from the author's work in geometry. The interest of this particular algebraic property is that it furnishes a sufficient, and presumably necessary, condition for another, geometric property of~$V$. Namely, the author obtained the following result:
\begin{theorem}{\cite{Smi16b}}
Let $G$~be a semisimple real Lie group, $V$ a representation of~$G$. Suppose that the action of~$w_0$ on~$V^L$ is nontrivial. Then there exists, in the affine group $G \ltimes V$, a subgroup $\Gamma$ whose linear part is Zariski-dense in~$G$, which is free of rank at least~$2$, and acts properly discontinuously on the affine space corresponding to~$V$.
\end{theorem}
He, and other people, also proved the converse statement in some special cases:
\begin{theorem}
The converse holds, for irreducible $V$:
\begin{itemize}
\item \cite{Smi18} if $G$~is split, but not of type $A_n$ ($n \geq 2$), $D_{2n+1}$ or~$E_6$;
\item \cite{Smi18} if $G$~is split, has one of these types, and $V$ satisfies a very restrictive additional assumption (see \cite{Smi18} for the precise statement);
\item \cite{AMS11} if $G = \SO(p,q)$ for arbitrary $p$ and~$q$, and $V = \RR^{p+q}$ is the standard representation.
\end{itemize}
\end{theorem}
Moreover, it seems plausible that, by combining the approaches of~\cite{Smi18} and~\cite{AMS11}, we might prove the converse in all generality. This geometric property is related to the so-called Auslander conjecture \cite{Aus64}, which is an important conjecture that has stood for more than fifty years and generated an enormous amount of work: see \eg \cite{Mil77, Mar83, FG83, AMS02, DGKPre} and many many others. For the statement of the conjecture as well as a more comprehensive survey of past work on it, we refer to~\cite{DDGS}.

\subsection{Statement of main results}

We have run some numerical experiments that allow us to conjecture the answer to Problem~\ref{nontrivial_w0_action_simple}.(ii) in the case where $\lie{g}^{\mathbb{C}}$ is of type $B_r$, $D_r$ or exceptional. Indeed, our numerical experiments prove this conjecture in some particular cases (see theorem~\ref{orth_and_excep_existing_data} below).

In the case where $\lie{g}^{\mathbb{C}}$ is of type $A_r$ or~$C_r$, we have also made some computations in low rank. Unfortunately, the data we have is not sufficient to be able to predict the general pattern (and we do not have enough computational power to generate more); see also the final remark in~\cite{Smi20} for more details.

\begin{conjecture}
\label{orth_and_excep_conj}
Assume that $\lie{g}^{\mathbb{C}}$ is of type $B_r$ (for some $r \geq 1$), $D_r$ (for some $r \geq 3$), $E_6$, $E_7$, $E_8$, $F_4$ or~$G_2$. Let $\lambda$ be a dominant integral weight of~$\lie{g}^{\mathbb{C}}$, $V_\lambda$ the irreducible (complex) representation of $\lie{g}$ with highest weight $\lambda$. Then:
\begin{hypothenum}
\item If $\lambda$ is one of the weights listed in Table~\ref{orth_and_excep_table}, then $\restr{w_0}{V_\lambda^L} = \pm \Id$.
\item If $V_\lambda^L \neq 0$ (this can be looked up in \cite[Table 1]{Smi20b}) but $\lambda$ does not occur in Table~\ref{orth_and_excep_table}, then $\restr{w_0}{V_\lambda^L} \neq \pm \Id$.
\end{hypothenum}
\end{conjecture}

Here are the cases in which we have checked this conjecture. (In fact, in a few special cases where we judged it to be useful and found it feasible, we have actually gone a bit farther than the cutoff figures listed below; but these details would be too tedious to list.)

\begin{proposition}~
\label{orth_and_excep_existing_data}
\begin{itemize}
\item Conjecture~\ref{orth_and_excep_conj}.(ii) holds for all real forms of $B_r$ with $r \leq 7$, of $D_r$ with $r \leq 9$, and of all exceptional algebras.
\item Conjecture~\ref{orth_and_excep_conj}.(i) holds for all the algebras in the same list, for weights $\lambda = \sum_{i=1}^r c_i \varpi_i$ satisfying $c_i \leq 3 p_i$ for all coefficients $c_i$, where $p_i$ is the least positive integer such that $V_{p_i \varpi_i}^L \neq 0$. (In the case of $i \in \{2p, 2p+1\}$ for $\lie{g} = \lie{so}(p,q)$ with $p = \frac{p+q-2}{4}$, this definition technically makes $p_i$ infinite, but we take the convention $p_i = 1$ instead.)
\end{itemize}
\end{proposition}

The proof of this proposition relies on the additivity property \cite[Proposition~1.(iii)]{Smi20}, which reduces it to a finite number of computations; and an algorithm to compute the restriction of $w_0$ to $V^L$ that the author has recently developed and implemented in the LiE software \cite{LiE}. The details of that algorithm will be published in a subsequent paper.

\subsection{Acknowledgements}

The author is supported by the European Research Council (ERC) under the European Union Horizon 2020 research and innovation programme (ERC starting grant DiGGeS, grant agreement No. 715982).

\begin{table}[b]
\caption{\label{orth_and_excep_table}Values of $\lambda$ for which $\restr{w_0}{V_\lambda^L} = \pm \Id$, for various algebras $\lie{g}$. The fundamental weights $\varpi_i$ are numbered using the Bourbaki ordering \cite{BouGAL456}. The coefficients $k$, $l$ and $m$ range in the nonnegative integers. Note that the lists may contain duplicates.}
\addtocounter{table}{-1}
\begin{subtable}{\textwidth}
\caption{\label{Br_subtable} Real forms of $B_r = \lie{so}_{2r+1}(\mathbb{C})$ ($r \geq 1$). In $\lie{so}(p,q)$, we always assume $p \leq q$.}
\centering\bigskip
\begin{tabular}[t]{llll}
\toprule
The algebra $\mathfrak{g}$ & Weights $\lambda$ & Conditions on indices & \parbox{2.3cm}{Conditions on\\ coefficients} \\
\midrule
\multirow{4}{*}{\parbox{2.3cm}{$\lie{so}(p,q)$ \\ ${}\quad p \leq \frac{p+q}{4}$ \\ ${}\quad p+q$ odd}}
& \multirow{4}{*}{$\lambda = k \varpi_i + l \varpi_{2p}$}
  & $i = 1$ or $2p-1$ & any $k$, any $l$ \\ 
& & $2 = i = 2p-2$ & any $k$, any $l$ \\ 
& & $i = 2$ or $2p-2$ & $k \leq 2$, any $l$ \\ 
& & $2 < i < 2p-2\;\wedge\; 2|i$ & $k \leq 1$, any $l$ \\
\midrule
\multirow{6}{*}{\parbox{2.3cm}{$\lie{so}(p,q)$ \\ ${}\quad p = \frac{p+q+1}{4}$}}
& \multirow{6}{*}{$\lambda = k \varpi_i + l \varpi_{q-p}$}
  & $i = 1$ & any $k$; $l \leq 2$ \\ 
& & $2 = i = q-p-1$ & any $k$; $l \leq 2$ \\ 
& & $2 = i < q-p-1$ & $k \leq 2,\; l \leq 2$ \\ 
& & $2 < i < q-p-1\;\wedge\; 2|i$ & $k \leq 1,\; l \leq 2$ \\
& & $2 < i = q-p-1$ & $k \leq 2,\; l \leq 2$ \\ 
& & $i = q-p$ & any $k$, any $l$ \\ 
\midrule
\multirow{8}{*}{\parbox{2.3cm}{$\lie{so}(p,q)$ \\ ${}\quad p > \frac{p+q+1}{4}$ \\ ${}\quad p+q$ odd}}
& \multirow{3}{*}{$\lambda = k \varpi_i + l \varpi_{q-p}$}
  & $i = 1$ & any $k$; $l \leq 1$ \\ 
& & $2 = i < q-p$ & $k \leq 2,\; l \leq 1$ \\ 
& & $2 < i < q-p\;\wedge\; 2|i$ & $k \leq 1,\; l \leq 1$ \\ \cmidrule{2-4}
& \multirow{5}{*}{$\lambda = k \varpi_i$}
  & $2 = i = \frac{p+q-1}{2}$ & any $k$ \\ 
& & $q-p < i = 2$ & $k \leq 2$ \\ 
& & $q-p < i$ & $k \leq 1$ \\ 
& & $q-p+1 = i = \frac{p+q-1}{2}$ & $k \leq 4$ \\ 
& & $q-p < i = \frac{p+q-1}{2}$ & $k \leq 2$ \\ 
\bottomrule
\end{tabular}
\end{subtable}
\end{table}

\clearpage

\begin{table}
\ContinuedFloat
\begin{subtable}{\textwidth}
\caption{\label{Dr_subtable} Real forms of $D_r = \lie{so}_{2r}(\mathbb{C})$ ($r \geq 3$). In $\lie{so}(p,q)$, we always assume $p \leq q$; and we denote by $r := \frac{p+q}{2}$ the (complex) rank.}
\centering\bigskip
\begin{tabular}[t]{llll}
\toprule
The algebra $\mathfrak{g}$ & Weights $\lambda$ & Conditions on indices & \parbox{2.3cm}{Conditions on\\ coefficients} \\
\midrule
\parbox{2.3cm}{$\lie{so}(p,q)$ \\ ${}\quad p \leq \frac{p+q}{4}-1$ \\ ${}\quad p+q$ even}
& \multicolumn{3}{l}{same as for $p \leq \frac{p+q}{4}$ in the $B_r$ case} \\
\midrule
\parbox{2.3cm}{$\lie{so}(p,q)$ \\ ${}\quad p = \frac{p+q-2}{4}$}
& \multicolumn{3}{l}{same as for $p \leq \frac{p+q}{4}$ in the $B_r$ case, with $\varpi_{2p}$ replaced by $(\varpi_{2p} + \varpi_{2p+1})$} \\
\midrule
\multirow{2}{*}{\parbox{2.3cm}{$\lie{so}(p,q)$ \\ ${}\quad p = \frac{p+q}{4}$}}
& \multicolumn{3}{l}{same as for $p \leq \frac{p+q}{4}$ in the $B_r$ case} \\ \cmidrule{2-4}
& \multicolumn{3}{l}{same as for $p \leq \frac{p+q}{4}$ in the $B_r$ case, with $\varpi_{2p}$ replaced by $\varpi_{2p-1}$} \\
\midrule
\multirow{6}{*}{\parbox{3.3cm}{$\lie{so}(p,q)$ \\ ${}\quad p > \frac{p+q}{4}$ \\ ${}\quad p+q \equiv 0 \pmod{4}$}}
& \multirow{6}{*}{$\lambda = k \varpi_i$}
  & $i = 1$ & any $k$ \\
& & $i = 2$ & $k \leq 2$ \\
& & $2 < i < r-1\;\wedge\; 2|i$ & $k \leq 1$ \\
& & $i \in \{r-1, r\}\;\wedge\; r = 4$ & any $k$ \\
& & $i \in \{r-1, r\}\;\wedge\; p = \frac{p+q}{4} + 1$ & $k \leq 4$ \\
& & $i \in \{r-1, r\}$ & $k \leq 2$ \\
\midrule
\multirow{4}{*}{\parbox{3.3cm}{$\lie{so}(p,q)$ \\ ${}\quad p > \frac{p+q}{4}$ \\ ${}\quad p+q \equiv 2 \pmod{4}$}}
& \multirow{4}{*}{$\lambda = k \varpi_i$}
  & $i = 1$ & any $k$ \\
& & $2 = i < q-p-1$ & $k \leq 2$ \\
& & $2 < i < q-p-1\;\wedge\; 2|i$ & $k \leq 1$ \\
& & $i \in \{r-1, r\}\;\wedge\; i \leq q-p+1$ & any $k$ \\
\midrule
$\lie{so}^*(6)$ & $\lambda = k\varpi_1 + l\varpi_i$ & $i \in \{2, 3\}$ & any $k$, any $l$ \\
\midrule
\multicolumn{4}{l}{$\lie{so}^*(8)$: see $\lie{so}(6,2)$, to which it is isomorphic} \\
\midrule
$\lie{so}^*(10)$ & $\lambda = k\varpi_i$ & $i \in \{1, 4, 5\}$ & any $k$ \\
\midrule
\multirow{3}{*}{$\lie{so}^*(12)$}
& \multirow{3}{*}{$\lambda = k \varpi_i$}
  & $i \in \{1, 2, 6\}$ & any $k$ \\
& & $i = 4$ & $k \leq 1$ \\
& & $i = 5$ & $k \leq 2$ \\
\midrule
\parbox{3.3cm}{$\lie{so}^*(2r)$ \\ ${}\quad r > 5$, $r$ odd} & $\lambda = k\varpi_i$ & $i = 1$ & any $k$ \\
\midrule
\multirow{4}{*}{\parbox{3.3cm}{$\lie{so}^*(2r)$ \\ ${}\quad r > 6$, $r$ even}}
& \multirow{4}{*}{$\lambda = k \varpi_i$}
  & $i = 1$ & any $k$ \\
& & $i = 2$ & $k \leq 2$ \\
& & $2 < i < r-1\;\wedge\; 2|i$ & $k \leq 1$ \\
& & $i \in \{r-1, r\}$ & $k \leq 4$ \\
\bottomrule
\end{tabular}
\end{subtable}
\end{table}

\clearpage

\begin{table}[t]
\ContinuedFloat
\begin{subtable}{\textwidth}
\caption{\label{excep_subtable} Real forms of exceptional algebras.}
\centering\bigskip
\begin{tabular}[t]{llll}
\toprule
The algebra $\mathfrak{g}$ & Weights $\lambda$ & Conditions on indices & \parbox{2.3cm}{Conditions on\\ coefficients} \\
\midrule
$\mathrm{E~I}$, $\mathrm{E~II}$ & $\lambda = 0$ & & \\
\midrule
$\mathrm{E~III}$ & $\lambda = k \varpi_i$ & $i \in \{1,6\}$ & any $k$ \\
\midrule
$\mathrm{E~IV}$ & $\lambda = k \varpi_2 + l \varpi_i$ & $i \in \{1,3,5,6\}$ & any $k$, any $l$ \\
\midrule
\midrule
\multirow{2}{*}{$\mathrm{E~V}$, $\mathrm{E~VI}$}
& $\lambda = k \varpi_i$ & $i = 1$ & $k \leq 2$ \\
& $\lambda = k \varpi_i$ & $i \in \{6, 7\}$ & $k \leq 1$ \\
\midrule
\multirow{2}{*}{$\mathrm{E~VII}$}
& $\lambda = k \varpi_i$ & $i = 1$ & any $k$ \\
& $\lambda = k \varpi_i$ & $i \in \{6, 7\}$ & $k \leq 1$ \\
\midrule
\midrule
\multirow{2}{*}{$\mathrm{E~VIII}$, $\mathrm{E~IX}$}
& $\lambda = k \varpi_i$ & $i = 1$ & $k \leq 1$ \\
& $\lambda = k \varpi_i$ & $i = 8$ & $k \leq 2$ \\
\midrule
\midrule
$\mathrm{F~I}$ & $\lambda = k \varpi_i$ & $i \in \{1, 4\}$ & $k \leq 2$ \\
\midrule
$\mathrm{F~II}$ & $\lambda = k \varpi_1 + l \varpi_2 + m \varpi_i$ & $i \in \{3, 4\}$ & any $k$, any $l$, any $m$ \\
\midrule
\midrule
$\mathrm{G}$ & $\lambda = k \varpi_i$ & $i \in \{1, 2\}$ & $k \leq 2$ \\
\bottomrule
\end{tabular}
\end{subtable}
\end{table}

\bibliographystyle{alpha}
\bibliography{/home/ilia/Documents/Travaux_mathematiques/mybibliography.bib}

\newcommand{\inprep}{in preparation}  \newcommand{\arx}[1]{arXiv:#1}
  \newcommand{\forthcm}{to appear}  \newcommand{\subm}{submitted}
  \newcommand{\noop}[1]{}
\begin{thebibliography}{DDGS22}

\bibitem[AMS02]{AMS02}
H.~Abels, G.~A. Margulis, and G.~A. Soifer.
\newblock On the {Z}ariski closure of the linear part of a properly
  discontinuous group of affine transformations.
\newblock {\em J. Differential Geom.}, 60:315--344, 2002.

\bibitem[AMS11]{AMS11}
H.~Abels, G.~A. Margulis, and G.~A. Soifer.
\newblock The linear part of an affine group acting properly discontinuously
  and leaving a quadratic form invariant.
\newblock {\em Geom. Dedicata}, 153:1--46, 2011.

\bibitem[Aus64]{Aus64}
L.~Auslander.
\newblock The structure of complete locally affine manifolds.
\newblock {\em Topology}, 3:131--139, 1964.

\bibitem[Bou68]{BouGAL456}
N.~Bourbaki.
\newblock {\em {\'E}l{\'e}ments de Math{\'e}matique, Groupes et Alg{\`e}bres de
  {L}ie : chapitres 4, 5 et 6}.
\newblock Hermann, 1968.

\bibitem[DDGS22]{DDGS}
J.~Danciger, T.~A. Drumm, W.~M. Goldman, and I.~Smilga.
\newblock Proper actions of discrete groups of affine transformations.
\newblock In {\em Dynamics, Geometry, Number Theory: the Impact of {M}argulis
  on Modern Mathematics}. University of Chicago Press, \noop{2022}\forthcm
  \noop{22}.
\newblock \arx{2002.09520}.

\bibitem[DGK16]{DGKPre}
J.~Danciger, F.~Gu{\'e}ritaud, and F.~Kassel.
\newblock Margulis spacetimes via the arc complex.
\newblock {\em Invent. Math.}, 204(1):133--193, 2016.

\bibitem[FG83]{FG83}
D.~Fried and W.~M. Goldman.
\newblock Three-dimensional affine crystallographic groups.
\newblock {\em Adv. in Math.}, 47:1--49, 1983.

\bibitem[LFS18]{LFlSm}
B.~Le~Floch and I.~Smilga.
\newblock Action of {W}eyl group on zero-weight space.
\newblock {\em C. R. Math. Acad. Sci. Paris}, 356(8):852--858, 2018.

\bibitem[Mar83]{Mar83}
G.~A. Margulis.
\newblock Free properly discontinuous groups of affine transformations.
\newblock {\em Dokl. Akad. Nauk SSSR}, 272:785--788, 1983.

\bibitem[Mil77]{Mil77}
J.~Milnor.
\newblock On fundamental groups of complete affinely flat manifolds.
\newblock {\em Adv. in Math.}, 25:178--187, 1977.

\bibitem[Smi20a]{Smi20}
I.~Smilga.
\newblock Action of the restricted {W}eyl group on the {$L$}-invariant vectors
  of a representation.
\newblock In V.~Dobrev, editor, {\em Proceedings of the XIII International
  Workshop ``Lie Theory and Its Applications in Physics'' (Varna, Bulgaria,
  June 2019)}, volume 335 of {\em Springer Proceedings in Mathematics and
  Statistics}, pages 365--372. Springer Singapore, 2020.

\bibitem[Smi20b]{Smi18}
I.~Smilga.
\newblock Construction of {M}ilnorian representations.
\newblock {\em Geom. Dedicata}, 206:55--73, 2020.

\bibitem[Smi22a]{Smi16b}
I.~Smilga.
\newblock Proper affine actions: a sufficient criterion.
\newblock {\em Math. Ann.}, \noop{2022}\forthcm \noop{22}.
\newblock \arx{1612.08942}.

\bibitem[Smi22b]{Smi20b}
I.~Smilga.
\newblock Representations having vectors fixed by a {L}evi subgroup.
\newblock {\em J. of Algebra}, \noop{2022}\forthcm \noop{22}.
\newblock \arx{2002.10928}.

\bibitem[vLCL00]{LiE}
M.~A.~A. van Leeuwen, A.~M. Cohen, and B.~Lisser.
\newblock {\em {LiE}, A Package for {L}ie Group Computations}, 2000.
\newblock \url{http://wwwmathlabo.univ-poitiers.fr/~maavl/LiE/}.

\end{thebibliography}
\end{document}